\documentclass[12pt]{article}

\usepackage{amssymb}
\usepackage{amsmath,amssymb,amsthm, amscd}
\usepackage{indentfirst}

\theoremstyle{definition}

\title{\textbf{Blaschke Isoparametric Hypersurfaces
 in the Conformal Space ${\mathbb Q}^{n+1}_1$}, II }
 \author {{\bf Tongzhu Li$^1$, Changxiong Nie$^2$} \\
\small{1 Department of Mathematics, Beijing Institute of
Technology,} \\
\small{Beijing, China,100081, E-mail: litz@bit.edu.cn} \\
\small{2 Faculty of Mathematics and Statistics, Hubei Key Laboratory of Applied Mathematics,  }\\
\small{Hubei University, Wuhan, China, 430062, E-mail: chxnie@163.com
}}

\begin{document}

\maketitle \footnotetext [1]{T. Z. Li is   supported by the
grant No. 11571037  of NSFC; } \footnotetext [2]{C. X. Nie is
partially supported
by the grant No. 11571037  of NSFC.}

\begin{abstract}
   Let $x:  \mathbf{M}  \rightarrow  {\mathbb Q}^{n+1}_1$ be a regular  space-like
   hypersurface   in the conformal space ${\mathbb Q}^{n+1}_1$.
     We classify all those  hypersurfaces with parallel Blaschke tensor in the conformal  space
      up to the conformal equivalence.
\end{abstract}



\par\noindent
{\bf {\S} 1. Introduction.}
\par\medskip

  In [7] we have classified all regular space-like hypersurfaces with two distinct constant Blaschke eigenvalues
   in the conformal space ${\mathbb Q}^{n+1}_1$. In this paper, we shall focus on a special class of
   Blaschke isoparametric hypersurfaces, which are called hypersurfaces with parallel Blaschke tensor.
   We shall find that those hypersurfaces with parallel Blaschke tensor must be of 1, 2, or 3 distinct
    constant  Blaschke eigenvalues.

 First of all,
 we shall present some concrete space-like Blaschke
isoparametric hypersurfaces
in $ {\mathbb Q}^{n+1}_1$ with two distinct Blaschke eigenvalues.
 The details of calculation will occur in the proof of the main theorem of this paper
 later.

{\bf Example 1.1.} For some real number $ r > 0 $ and integers $n,k$ satisfying
$k=1,\cdots,n-1$, $u:   \mathbf N^k \rightarrow \mathbb S^{ k  +  1 }_1(r)\subset
\mathbb R^{ k  +  2 }_1$
 be a regular maximal space-like hypersurface  with constant scalar curvature
$$ \rho_1 = \frac{  k ( k - 1 ) }{  r^2 }  -  \frac{  n - 1 }{ n }. $$
Let
 $\{e_1, \cdots , e_k\}$ be an local basis for $u$ with dual basis
 $\{\omega^1, \cdots , \omega^k\}$. Denote the second
fundamental form of $ u $ by $II _1 =  \sum_ { ij } h _{ ij }
\omega^i\otimes \omega^j$. Then $$ \sum_{ i } h^i_i = 0,\quad
\sum_ { ij } h^i_j h^j_i = \frac{  n - 1 }{ n },$$
 where the  indices are shifted by the first fundamental form $ I_1 $.
Denote
 $\Delta_{ 1 }$ the Laplacian   with respect to $I_1$. It is easy to know that
 $$\Delta_{ 1 } u  =  -  \frac{  k u }{ r^2 }.$$
 Let $ v:  \mathbb H^{ n  -  k  } ( r ) \rightarrow \mathbb R^{ n  -  k  +  1 }_1$
 be the standard totally umbilical hypersurface.
 Then the scalar curvature
$$ \rho_2 = - \frac{  ( n - k  )  ( n - k - 1 )  }{  r^2 } . $$
 And we have
$$\Delta_{ 2 } v  =  \frac{  ( n - k  ) v }{ r^2 },$$
where  $\Delta_{ 2 }$ is the Laplacian   with respect to the first fundamental form  $ I_2 $
 of $v$.
For  $  y = ( u , v ) :   \mathbf M  =  \mathbf N^k \times \mathbb H^{ n  -  k  } ( r )
\rightarrow C^{ n + 2 } \subset  \mathbb R^{n+ 3 }_2$, we find exactly that the conformal metric
 $ g = \langle \mathrm d y , \mathrm d y \rangle$. Therefore $y$ is the canonical
  lift of $x = [ y ] :  \mathbf{M}  \rightarrow  {\mathbb Q}^{n+1}_1$.
  A direct calculation will yield that
  $$ ( A^i_j ) = \frac{ 1 }{ 2r^2 } (\mathrm I_k \oplus ( - \mathrm I_{ n - k } ) ) ,\quad
  ( B^i_j ) = ( h^i_j ) \oplus   \mathbf 0_{ n - k }  ,\quad
   C_i=0,\forall i,$$
 where $\mathrm I_k$ means $k$-rank
 unit matrix and $ \mathbf 0_k$ means $k$-rank
 zero matrix.

{\bf Example 1.2.} For some real number $ r > 0 $ and integers $n,k$ satisfying
$k=1,\cdots,n-1$, $u:   \mathbf N^k \rightarrow \mathbb H^{ k  +  1 }_1(r)\subset
\mathbb R^{ k  +  2 }_2$
 be a regular maximal space-like hypersurface  with constant scalar curvature
$$ \rho_1 = - \frac{  k ( k - 1 ) }{  r^2 }  -  \frac{  n - 1 }{ n }. $$
Let
 $\{e_1, \cdots , e_k\}$ be an local basis for $u$ with dual basis
 $\{\omega^1, \cdots , \omega^k\}$. Denote the second
fundamental form of $ u $ by $II _1 =  \sum_ { ij } h _{ ij }
\omega^i\otimes \omega^j$. Then $$ \sum_{ i } h^i_i = 0,\quad
\sum_ { ij } h^i_j h^j_i = \frac{  n - 1 }{ n }.$$

Denote
 $\Delta_{ 1 }$ the Laplacian   for $I_1$. It is easy to know that
 $$\Delta_{ 1 } u  =  \frac{  k u }{ r^2 }.$$
 Let $ v:  \mathbb S^{ n  -  k  } ( r ) \rightarrow \mathbb R^{ n  -  k  +  1 }$
 be the standard totally umbilical hypersurface.
 Then the scalar curvature
$$ \rho_2 = \frac{  ( n - k  )  ( n - k - 1 )  }{  r^2 } . $$
 And we have
$$\Delta_{ 2 } v  =  - \frac{  ( n - k  ) v }{ r^2 },$$
where  $\Delta_{ 2 }$ is the Laplacian   with respect to the first fundamental form  $ I_2 $
 of $v$.
For  $  y = ( u , v ) :   \mathbf M  =  \mathbf N^k \times \mathbb S^{ n  -  k  } ( r )
\rightarrow C^{ n + 2 } \subset  \mathbb R^{n+ 3 }_2$,  we find exactly that the conformal metric
 $ g = \langle \mathrm d y , \mathrm d y \rangle$. Therefore $y$ is the canonical
  lift of $x = [ y ] :  \mathbf{M}  \rightarrow  {\mathbb Q}^{n+1}_1$.
  A direct calculation  will yield that
  $$ ( A^i_j ) = - \frac{ 1 }{2 r^2 } ( I_k \oplus ( - I_{ n - k } ) ) ,\quad
  ( B^i_j ) = ( h^i_j ) \oplus ( \mathbf 0_{ n - k } ),\quad
   C_i=0,\forall i.$$

Now the
main theorem in this paper is stated as follows:

{\bf Theorem B} Let $x:  \mathbf{M}  \rightarrow  {\mathbb Q}^{n+1}_1$ be a regular
 space-like  hypersurface  in the conformal space ${\mathbb Q}^{n+1}_1$.
If the Blaschke tensor $\mathbb{A}$ of $x$ is parallel,
 then one of
the following statements holds:

(1) $x$ is conformal isotropic and is therefore locally conformally equivalent to:

\quad (a) a maximal hypersurface in ${\mathbb S}^{n+1}_1$ with constant scalar curvature; or

\quad (b)   a maximal hypersurface in ${\mathbb R}^{n+1}_1$ with constant
scalar curvature; or

\quad (c)   a maximal hypersurface in ${\mathbb H}^{n+1}_1$ with
constant scalar curvature;

(2) $x$ is of parallel conformal second fundamental form $\mathbb B$ and is therefore locally
conformally
equivalent to:

\quad  (a) a standard cylinder ${\mathbb H}^{k} (r)
\times{\mathbb S}^{n-k}(\sqrt{1+r^2}),r>0$, in ${\mathbb S}^{n+1}_1$ for some  positive
integer $k$ and $n-k$; or

\quad (b)  a standard cylinder   ${\mathbb H}^{k}
\times{\mathbb R}^{n-k} $ in ${\mathbb R}^{n+1}_1$ for some   positive
integer $k$ and $n-k$; or

\quad (c)   a standard cylinder ${\mathbb H}^{k} (r)
\times{\mathbb H}^{n-k}(\sqrt{1-r^2}),0<r<1$,
 in ${\mathbb H}^{n+1}_1$ for some   positive
integer $k$ and $n-k$; or

\quad (d)
 the wrapped product embedding
 $$
 u:\mathbb H^{p}(r)\times
 \mathbb S^{q}(\sqrt{ r^2+1})\times
 \mathbb R^+\times\mathbb R^{n-p-q-1}\subset\mathbb R^{n+2}_1 \rightarrow
 \mathbb R^{n+1}_1, $$
 $$
  (u',u'',t,u''')\mapsto(tu',tu'',u'''),
 $$
  where
  $$u'\in \mathbb H^{p}(r),u''\in\mathbb S^{q}(\sqrt{ r^2+1}),t>0,
  u'''\in\mathbb R^{n-p-q-1},r>0,$$
for some   positive
integer $p$, $q$, and $n-p-q-1$;

(3) $x$ is non-conformal isotropic with a non-parallel conformal second fundamental form B and is
locally conformally equivalent to:

\quad (a) one of the  hypersurfaces as indicated in Example 1.1; or

 \quad (b) one of the   hypersurfaces as indicated in Example 1.2.

{\bf Remark 1.4}. All the  above examples are Blaschke isoparametric. Among them, case (1) is of one eigenvalue of the
 Blaschke tensor $\mathbb A$, the first three examples of case (2) and case (3) are of two eigenvalues,
 and subcase (d) of case (2) is of three eigenvalues.

This paper is organized as follows. In Section 2 we will
 give some frequently-used equations occurred many times in [8-10] {\it etc.}, especially in [7]. Readers may find more details
 in those conferences.  In additional, we prove some lemmas which will be used in the next section.
 In
Section 3  we  prove the Theorem B.

\bigskip

\medskip
\par\noindent
{\bf {\S} 2. The fundamental equations and some lemmas.}
\par\medskip

Let $\mathbb R^{N}_{s}$ denote pseudo-Euclidean space, which is the real vector space $\mathbb R^{N}$
with the non-degenerate inner product $\langle,\rangle$ given by
  $$\langle\xi,\eta\rangle
  = - \sum_{i=1}^{ s } x_iy_i  +  \sum_{ i =  s + 1 }^{N} x_iy_i,
  \eqno{(1. 1)}$$
  where $\xi=(x_1, \cdots
x_{ _{N} } ),\eta=(y_1, \cdots , y_{ _{N} } )\in\mathbb R^{N}$.

Let
$$C^{n+1}:=\{\xi\in\mathbb R^{n+2}_{s+1}|\langle \xi  ,\xi  \rangle=0,\xi  \neq0\}, \eqno{(1.2)}$$
 $$\mathbb Q^n_s:=\{[\xi  ]\in\mathbb R P^{n+1}|\langle \xi  ,\xi  \rangle=0\}=
 C^{n+1}/(\mathbb R\backslash\{0\}).\eqno{(1.3)}$$
 We call
$C^{n+1}$ the light cone in $\mathbb R^{n+2}_{s+1}$ and  $\mathbb Q^n_s$ the
conformal space ({\sl or} projective light cone) in $\mathbb R P^{n+1}$.

  The standard metric $h$ of the conformal space $\mathbb Q^n_s$ can be obtained
 through the pseudo-Riemannian submersion
 $$\pi:C^{n+1}\rightarrow\mathbb Q^n_s,\xi\mapsto[\xi].$$
   We can check $(\mathbb Q^n_s,h)$ is a pseudo-Riemannian manifold.

  We define the pseudo-Riemannian sphere space $\mathbb S^n_s ( r ) $  and
pseudo-Riemannian hyperbolic space  $\mathbb H^n_s ( r ) $ with radius $r$ by
$$\mathbb S^n_s ( r ) =\{u\in \mathbb R^{n+1}_s|\langle u, u\rangle = r^2\},  \quad
\mathbb H^n_s ( r ) =\{u\in \mathbb R^{n+1}_{s+1}|\langle u, u\rangle = - r^2\}.   $$
When   $ r =1 $ we usually omit the radius $ r   $. When $s=1$ and $ r =1 $ we call them de Sitter space
$\mathbb S^n_1$ and anti-de Sitter space $\mathbb H^n_1$.

We may assume    ${\mathbb Q}^n_s$ as the common compactification of
$\mathbb R^n_s$, $\mathbb S^n_s$ and $ \mathbb H^n_s$, and
 $\mathbb R^n_s$, $\mathbb S^n_s$ and $\mathbb H^n_s$ as the subsets
 of $\mathbb
 Q^n_s$ when referring to  the conformal geometry. 

Let $x:  \mathbf{M}  \rightarrow  {\mathbb Q}^{n+1}_1$ be a regular
  space-like hypersurface  in the conformal space ${\mathbb Q}^{n+1}_1$.
We skip the standard
 procedure of the analysis of the conformal geometry of hypersurfaces.
 For more details, readers may see [8], {\it etc.}. We have four important conformal invariants,
 the conformal metric $g$,
 the Blaschke tensor $\mathbb A$,
 the conformal second fundamental form
  $\mathbb
B$,
 the conformal   form $\phi
 $. Then we have some   fundamental equations  of $x$ that are used later as follows:
  $$A_{ij, k}-A_{ik, j} = B_{ij}C _k-B_{ik}C_j, \eqno{(2.1)}$$
  $$ B_{ij,
  k}-B_{ik, j}=  g_{_{ij}}C_k - g_{_{ik}}C_j,
  \eqno{(2.2)}$$
  $$C_{i, j}-C_{j, i}
  =\sum_{kl}g^{kl}(B_{ik}A_{lj}-B_{jk}A_{li}),
   \eqno{(2.3)}$$
  $$R_{ijkl}=
    (g_{ik}A_{jl}-g_{il}A_{jk})
    +(A_{ik}g_{jl}-A_{il}g_{jk})
    -(B _{ik}B_{jl}-B _{il}B _{jk}),\eqno{(2.4)}$$
  $$\sum_{ij } B^i_jB^ j_i=     \frac{n-1}{n},\eqno{(2.5)}$$
  $$
\sum_i B^i_i=0. \eqno{(2.6)}$$

  For  $u:\mathbf
M\rightarrow \mathbf L^{n+1}(\epsilon)$,
when $\epsilon=0,1,-1,\mathbf L^{n+1}(\epsilon)$ denotes $\mathbf
R^{n+1}_1$, $\mathbf
S^{n+1}_1$ and $\mathbf
H^{n+1}_1$, respectively. We have
 $$
  e^{2\tau} = \frac{1 }{n-1}(  n  \sum_{ij}h^i_jh^j_i -
  ( \sum_i h^i_i  ) ^2) .
  \eqno{(2.7)}$$
  $$    A_{ij} = \tau_i\tau_j -
 H h_{ij}  -\tau_{i,
  j}  -  \frac{1}{2}(\sum_i \tau^i \tau_i-H^2 - \epsilon ) I_{ij},  \eqno{(2.8)}$$
  $$  B_{ij}  = e^{ \tau} ( h_{ij}-H I_{ij} ) ,
  \eqno{(2.9)}$$
  $$ C_i=   e^{- \tau}  ( H\tau_i
 -\sum_{j}h_{ij}\tau^j-H_{ i} )
 .\eqno{(2.10)}$$

  We remind readers that if we call a tensor is parallel
  then usually we regard the first order covariant differential of
  the tensor vanishes. So the Blaschke tensor $\mathbb A$ of the regular space-like hypersurface $x$
   is defined parallel if and only if $A_{ij,k}=0$,
   for any $i,j,k=1,\cdots,n$.

  Next we introduce several lemmas.

 {\bf Lemma 2.1.}\ Let $x:  \mathbf{M}  \rightarrow  {\mathbb Q}^{n+1}_1$ be a regular
 space-like   hypersurface.
 If the Blaschke tensor $\mathbb A$ is parallel, then the conformal form $\phi$ vanishes.

{\sl Proof}\quad
  Denote $B = ( B_{ij}) $.
  If we choose an appropriate orthonormal basis $ e_1, \cdots , e_n $,
  we can write
  $$ B = ( b_1 I_{ m_1 } ) \oplus\cdots\oplus ( b_s I_{ m_s } ) ,$$
  where $b_1,\cdots,b_s$ are  distinct conformal eigenvalues
  and $m_1,\cdots,m_s$ are some positive integers satisfying $ \sum_{t=1}^sm_t=n$.

 From now on  we adopt the convention on the ranges
of indices in this section:
  $$1\leq i, j, k, l\leq n,\quad 1\leq t,t'\leq s
  . $$
  Denote $$\Gamma_t=\{i|B_{ii}=b_t\}.$$
  From the equation (2.1) we have
  $$B_{ij}C_k=B_{ik}C_j.\eqno{(2.11)}$$
 By (2.5) and (2.6), there must be
  at least two distinct $t,t'$ such that $b_t,b_{t'}\neq0$.
  Taking some fixed $i=j\in\Gamma_t$ in (2.11), we shall find that all the $C_k$'s are zero  except $C_i$.
  On the other hand, taking some fixed $i'=j'\in\Gamma_{t'}$ in (2.11), we shall find that all the $C_k$'s are zero  except $C_{i'}$.
  From above we know that
  all the $C_k$'s are zero. That means, the conformal form $\phi$ vanishes. $\Box$

  {\bf Lemma 2.2.}\ Let $x:  \mathbf{M}  \rightarrow  {\mathbb Q}^{n+1}_1$ be a regular
 space-like Blaschke isoparametric hypersurface.
 If the Blaschke tensor $\mathbb A$ is parallel, then $x$ is Blaschke isoparametric.
 And    the number
 of distinct Blaschke
 eigenvalues is no greater than 3.

{\sl Proof}\quad
  Denote $ A = ( A_{ij} ) $.
  If we choose an appropriate orthonormal basis $ e_1, \cdots , e_n $,
  we can write
  $$A = ( a_1 I_{ m_1 } ) \oplus\cdots\oplus ( a_s I_{ m_s } ) ,$$
  where $a_1,\cdots,a_s$ are  distinct Blaschke eigenvalues
  and $m_1,\cdots,m_s$ are some positive integers satisfying $ \sum_{t=1}^sm_t=n$.

 From now on  we adopt the convention on the ranges
of indices in this section:
  $$1\leq i, j, k, l\leq n,\quad 1\leq t,t',t''\leq s
  . $$
  Denote $$\Gamma_t=\{i|A_{ii}=a_t\}.$$
  Taking  $i=j\in\Gamma_t$ in
  $$  \sum_k A_{ij,k}\omega^k=\mathrm d A_{ij}  -  \sum_k A_{kj}\omega^k_i
   -  \sum_k A_{ik}\omega^k_j,\eqno{(2.12)}$$
   we shall get $\mathrm da_t=0$ because of $\omega^i_j+\omega^ j_i=0$, which implies that $a_t$'s
   are constant. Therefore  $x$ is Blaschke isoparametric.

   From the Lemma 2.1, we know that the conformal form $\phi=0$.
   Therefore by  (2.3), we have  $[A,B]=0$.
   From the algebraic and geometric technics, we can modify the orthonormal basis $ e_1, \cdots , e_n $
 such that the matrix $B$ can be diagonalized into the form
 $ B = \mathrm{diag}( b_1,\cdots,b_n )$.

  Taking  $i \in\Gamma_t,j\in\Gamma_{t'},t\neq t'$ in (2.12), we get
  $$\omega^i_j=0,\ \mathrm{when}\ i \in\Gamma_t,j\in\Gamma_{t'},t\neq t'.\eqno{(2.13)}$$
  The above equation will yield
  $$0=\mathrm d \omega^i_j + \sum_k \omega^i_k\wedge\omega^k_j=\Omega^i_j,$$
  which implies by (2.4) that
  $$R_{ijij}=   a_t + a_{t'} - b_i b_j =0
  ,\ \mathrm{when}\ i \in\Gamma_t,j\in\Gamma_{t'},t\neq t'. \eqno{(2.14)}$$
  Fixing one  of  two subscripts $i$ and $j$ in (2.14) and letting the other subscript vary,
  we will easily obtain $ B = ( b_1 I_{ m_1 } ) \oplus\cdots\oplus ( b_s I_{ m_s } ) $.
  So (2.14) comes to
  $$    a_t + a_{t'} - b_t b_{t'} =0
  ,\ \mathrm{where}\  t\neq t'. \eqno{(2.15)}$$

  If the number
 of distinct Blaschke
 eigenvalues $s>3$, it is easy to induce by (2.15) that,
 for distinct three number $t, t',t''$,
 $$ a_t-a_{t'} =b_{t''}  ( b_t-b_{t'} ), \eqno{(2.16)}$$
 $$\frac{a_t-a_{t'}}{b_t-b_{t'}}=b_{t''}. \eqno{(2.17)}$$
  In fact, we can guarantee that when $t\neq t'$, $b_t\neq b_{t'}$ by (2.16).
  From (2.17) we know that for fixed $t,t'$, the number
  $\frac{a_t-a_{t'}}{b_t-b_{t'}}$ will have more than 2 values.
  That is a contraction. So the assumption $s>3$ is wrong.  The number
 of distinct Blaschke
 eigenvalues is no greater than 3. $\Box$

  {\bf Lemma 2.3.}\ Let $x:  \mathbf{M}  \rightarrow  {\mathbb Q}^{n+1}_1$ be a regular
 space-like Blaschke isoparametric hypersurface with parallel Blaschke tensor.
 If $x$ has three distinct Blaschke
 eigenvalues, then the conformal second fundamental form $\mathbb B$ is parallel.

{\sl Proof}\quad
   From the proof of the Lemma 2.2,  we can  choose an appropriate orthonormal basis $ e_1, \cdots , e_n $
   such that
  $$A = ( a_1 I_{ m_1 } ) \oplus( a_2 I_{ m_2 } )\oplus ( a_3 I_{ m_3 } ) ,
  B= (b_1 I_{ m_1 } ) \oplus( b_2 I_{ m_2 } )\oplus ( b_3 I_{ m_3 } ) ,$$
  where $a_1,a_2,a_3$ are  distinct constant Blaschke eigenvalues.

  By (2.16) we know that $b_1,b_2, b_3$ are distinct and non-zero.
  Combining (2.5), (2.6) and (2.15), we can compute that $b_1,b_2, b_3$ are constant.

   For $t=1,2,3$, taking  $i,j\in\Gamma_t$ in
  $$ \sum_k B_{ij,k}\omega^k=\mathrm d B_{ij}  - \sum_k B_{kj}\omega^k_i
   -  \sum_k B_{ik}\omega^k_j,\eqno{(2.18)}$$
    we shall get
   $$B_{ij,k}=0,\ \mathrm{ where}\ i,j\in\Gamma_t,\forall t, k.\eqno{(2.19)}$$
  Taking $i\in\Gamma_1,j\in\Gamma_2$ in (2.18), and recalling (2.13),
  we obtain
   $$B_{ij,k}=0,\ \mathrm{ where}\ i\in\Gamma_1,j\in\Gamma_2, k\in\Gamma_3.\eqno{(2.20)}$$

  From Lemma 2.1 and (2.2), we know that $B_{ij,k}$ is all symmetric with respect to the subscripts.
  Combining (2.19) and (2.20), we know that $B_{ij,k}=0, \forall i,j, k. $
  $\Box$

\bigskip

\par\noindent
{\bf {\S} 3. Proof of the {\bf Theorem B}.}
\par\medskip
 Let $x:  \mathbf{M}  \rightarrow  {\mathbb Q}^{n+1}_1$ be a regular
 space-like Blaschke isoparametric hypersurface with parallel Blaschke tensor.
 From the Lemma 2.2, we know that  the number
 of distinct Blaschke
 eigenvalues is no greater than 3.

  It suffices to consider the following three cases:

Case (I): the number
 of distinct Blaschke
 eigenvalues is 1.

 In this case, we know that $x$ is actually  conformal
 isotropic. In the Theorem 5.2 of [8], we have classified the conformal
 isotropic submanifolds in the conformal space ${\mathbb Q}^n_p$. So we have the

 {\bf Theorem 3.1}\ Let $x:  \mathbf{M}  \rightarrow  {\mathbb Q}^{n+1}_1$ be a regular
 space-like  hypersurface  in the conformal space ${\mathbb Q}^{n+1}_1$.
If   $x$ is conformal
 isotropic,
 then   $x$ is   locally conformally equivalent to:

\quad (a) a maximal hypersurface in ${\mathbb S}^{n+1}_1$ with constant scalar curvature; or

\quad (b)   a maximal hypersurface in ${\mathbb R}^{n+1}_1$ with constant
scalar curvature; or

\quad (c)   a maximal hypersurface in ${\mathbb H}^{n+1}_1$ with
constant scalar curvature.

Case (II):   the number
 of distinct Blaschke
 eigenvalues is 2.

 In this case, we note that   we have classified all the space-like Blaschke
isoparametric hypersurfaces with two distinct Blaschke eigenvalues
 in the Theorem A of [11]. For the matter of completeness we state the following  modifying version
 of    the Theorem A. We have the

 {\bf Theorem 3.2}\ Let $x:  \mathbf{M}  \rightarrow  {\mathbb Q}^{n+1}_1$ be a regular
 space-like  hypersurface  in the conformal space ${\mathbb Q}^{n+1}_1$.
If $x$ is of two distinct constant Blaschke eigenvalues and of vanishing conformal form, then one of
the following statements holds:

(1) $x$ is conformal isoparametric with two distinct conformal principal  curvatures and is therefore
locally conformally equivalent to the regular
 space-like  hypersurface  in   ${\mathbb Q}^{n+1}_1$ determined by:

\quad  (a) a standard cylinder ${\mathbb H}^{k} (r)
\times{\mathbb S}^{n-k}(\sqrt{1+r^2}),r>0$ in ${\mathbb S}^{n+1}_1$ for some  positive
integer $k$ and $n-k$; or

\quad (b)  a standard cylinder   ${\mathbb H}^{k}
\times{\mathbb R}^{n-k} $ in ${\mathbb R}^{n+1}_1$ for some   positive
integer $k$ and $n-k$; or

\quad (c)   a standard cylinder ${\mathbb H}^{k} (r)
\times{\mathbb H}^{n-k}(\sqrt{1-r^2}),0<r<1$,
 in ${\mathbb H}^{n+1}_1$ for some   positive
integer $k$ and $n-k$;

(2) $x$   is locally conformally equivalent to:

\quad (a) one of the   hypersurfaces as indicated in  Example 1.1;  or

\quad (b) one of the   hypersurfaces as indicated in  Example 1.2.

Case (III):   the number
 of distinct Blaschke
 eigenvalues is 3.

 In this case,  from the lemma 2.3 we know that the conformal second fundamental form $\mathbb B$ is parallel.
 But we have classified all the space-like hypersurfaces with parallel conformal second
fundamental forms
 in the Classification Theorem   of [10]. For the matter of completeness we state the following  modifying version
 of    the Classification Theorem. We have the

 {\bf Theorem 3.3}\ Let $x:  \mathbf{M}  \rightarrow  {\mathbb Q}^{n+1}_1$ be a regular
 space-like  hypersurface  in the conformal space ${\mathbb Q}^{n+1}_1$. If
 $x$ is of parallel conformal second fundamental form $\mathbb B$, then it is locally
conformally
equivalent to:

\quad  (a) a standard cylinder ${\mathbb H}^{k} (r)
\times{\mathbb S}^{n-k}(\sqrt{1+r^2}),r>0$, in ${\mathbb S}^{n+1}_1$ for some  positive
integer $k$ and $n-k$; or

\quad (b)  a standard cylinder   ${\mathbb H}^{k}
\times{\mathbb R}^{n-k} $ in ${\mathbb R}^{n+1}_1$ for some   positive
integer $k$ and $n-k$; or

\quad (c)   a standard cylinder ${\mathbb H}^{k} (r)
\times{\mathbb H}^{n-k}(\sqrt{1-r^2}),0<r<1$,
 in ${\mathbb H}^{n+1}_1$ for some   positive
integer $k$ and $n-k$; or

\quad (d)
 the wrapped product embedding
 $$
 u:\mathbb H^{p}(r)\times
 \mathbb S^{q}(\sqrt{ r^2+1})\times
 \mathbb R^+\times\mathbb R^{n-p-q-1}\subset\mathbb R^{n+2}_1 \rightarrow
 \mathbb R^{n+1}_1, $$
 $$
  (u',u'',t,u''')\mapsto(tu',tu'',u'''),
 $$
  where
  $$u'\in \mathbb H^{p}(r),u''\in\mathbb S^{q}(\sqrt{ r^2+1}),t>0,
  u'''\in\mathbb R^{n-p-q-1},r>0,$$
for some   positive
integer $p$, $q$, and $n-p-q-1$.

   The last thing we ought to prove is that
 the wrapped product embedding $u$ has 3 distinct constant Blaschke
 eigenvalues. In fact,
 the first fundamental form of $u$ is
 $$
 I=\langle\mbox d{}u,\mbox d{}u\rangle =
 t^2\mbox d{}u'\cdot\mbox d{}u'
 +
 t^2\langle\mbox d{}u'',\mbox d{}u''\rangle+\mbox d{}t\otimes\mbox d{}t+\mbox d{}u'''\cdot\mbox d{}u'''.$$
 The unit time-like normal vector field of $u$
 $$
 e_{n+1}=(\frac{r}{\sqrt{r^2+1}}u',\frac{\sqrt{r^2+1}}{r}u'',\mathbf 0).$$
 The second fundamental form of $u$ is
 $$
 II=\langle\mbox d{}u,\mbox d{}e_{n+1}\rangle=
 t( \frac{r}{\sqrt{r^2+1}}\mbox d{}u'\cdot\mbox d{}u'+\frac{\sqrt{r^2+1}}{r}
 \langle\mbox d{}u'' ,\mbox d{}u'' \rangle ).$$
 Denote $I_k$ as $k$ order unit matrix, $\mathbf{0}_k$ as $k$ order zero matrix.
 Then we have $$
 (h^i_{ j})=( \frac{r}{\sqrt{r^2+1}t}I_p)\oplus
 ( \frac{\sqrt{r^2+1}}{rt}I_q)\oplus
 \mathbf{0}_{n-p-q},\quad
 H= \frac{pr^2+q(r^2+1)}{nr\sqrt{r^2+1}t},\eqno( 3.1)$$
 $$
 e^{2\tau} =\frac{p(n-p)r^4-2pqr^2(r^2+1)+q(n-q)(r^2+1)^2}{n-1}\cdot\frac{1}{t^2}
 :=\frac{d}{t^2}.\eqno(3.2)$$
 A directive calculation tells us that
 $$\tau_{i}=0,  i\neq p+q+1;\quad
 \tau_{p+q+1}=-\frac{1}{t },\eqno(3.3)$$
 $$\tau_{i,j}=0,(i,j)\neq (p+q+1,p+q+1);\quad \tau_{p+q+1,p+q+1}=\frac{1}{t^2}.\eqno(3.4)$$
 From (2.8) we get
 $$A^i_j
 =\sum_k e^{-2\tau} I^{ik} A_{kj}
 =      e^{-2\tau}     (   \tau^i\tau_j -
 H h^i_j  -\tau^i_{,j} -  \frac{1}{2}(\sum_k \tau^k \tau_k-H^2   )\delta^i_j).\eqno(3.5)$$
 Let
 $$a= \frac{r}{\sqrt{r^2+1} },\quad
 b=\frac{\sqrt{r^2+1}}{r},\quad
 c=\frac{pr^2+q(r^2+1)}{nr\sqrt{r^2+1}},$$
 $$
 d=\frac{p(n-p)r^4-2pqr^2(r^2+1)+q(n-q)(r^2+1)^2}{n-1}.$$
 From (3.1)-(3.5) we know that $u$ has three constant Blaschke eigenvalues
 $$a_1=\frac{c^2-2a-1}{2d},\quad a_2=\frac{c^2-2b-1}{2d},\quad a_3=\frac{c^2 -1}{2d}.$$
  Therefore we have the

 {\bf Theorem 3.4}\ Let $x:  \mathbf{M}  \rightarrow  {\mathbb Q}^{n+1}_1$ be a regular
 space-like  hypersurface  in the conformal space ${\mathbb Q}^{n+1}_1$. If
 $x$ is of parallel conformal second fundamental form $\mathbb B$ and of three constant Blaschke eigenvalues, then it is locally
conformally
equivalent to
 the wrapped product embedding
 $$
 u:\mathbb H^{p}(r)\times
 \mathbb S^{q}(\sqrt{ r^2+1})\times
 \mathbb R^+\times\mathbb R^{n-p-q-1}\subset\mathbb R^{n+2}_1 \rightarrow
 \mathbb R^{n+1}_1, $$
 $$
  (u',u'',t,u''')\mapsto(tu',tu'',u'''),
 $$
  where
  $$u'\in \mathbb H^{p}(r),u''\in\mathbb S^{q}(\sqrt{ r^2+1}),t>0,
  u'''\in\mathbb R^{n-p-q-1},r>0,$$
for some   positive
integer $p$, $q$, and $n-p-q-1$.

 Summing up the above process, especially using the Theorems 3.1, 3.2, 3.4,
 we have proved
 the Theorem B.


\end{document}